\documentclass[11pt]{amsart} 
\usepackage{amsmath}
\usepackage{amsthm}
\usepackage{amssymb}
\usepackage{latexsym}
\usepackage{chemarr}\usepackage{setspace}
\usepackage{mathrsfs}
\usepackage{bbm}
\usepackage[usenames]{color}

\newtheorem{definition}{Definition}
\newtheorem{theorem}{Theorem}
\newtheorem{proposition}[theorem]{Proposition}
\newtheorem{lemma}{Lemma}[theorem]

\newcommand{\ra}{\rightarrow}
\def\epsilon{\varepsilon}

\def \ul {\underline}
\def \bl {\begin{lemma}}
\def \el {\end{lemma}}
\def \bt {\begin{theorem}}
\def \et {\end{theorem}} 

\def \Tkk {\mathcal T _{k+1}}
\def \mc{\mathcal}
\def \limn {\lim_{n \ra \infty}}
\def \limk {\lim_{k \ra \infty}}

\def \Tk{\mathcal T _k}
\def \Ck{\mathcal C _k}
\def \IR{\mathbb R}
\def \IN{\mathbb N}

\def \htop{h_{top}}
\def \bp{\begin{proof}}
\def \ep{\end{proof}}
\def \be {\begin{equation}}
\def \ee {\end{equation}}
\def \beq {\begin{eqnarray*}}
\def \eeq {\end{eqnarray*}}

\def \line{|}

\def \bd {\begin{definition}}
\def \ed {\end{definition}}
\def \xr {X_\rho}
 
\begin{document}

\title[A Variational Principle for Pressure]{A Variational Principle for Topological Pressure for certain non-compact sets}
\author{Dan Thompson, University of Warwick}
\date{23rd September 2008}
\begin{abstract}
Let $(X,d)$ be a compact metric space, $f:X \mapsto X$ be a continuous map with the specification property, and $\varphi: X \mapsto \IR$ be a continuous function. We prove a variational principle for topological pressure (in the sense of Pesin and Pitskel) for non-compact sets of the form
\[
\left \{ x \in X : \lim_{n \ra \infty} \frac{1}{n} \sum_{i = 0}^{n-1} \varphi (f^i (x)) = \alpha \right \}.
\]
Analogous results were previously known for topological entropy. As an application, we prove multifractal analysis results for the entropy spectrum of a suspension flow over a continuous map with specification and the dimension spectrum of certain non-uniformly expanding interval maps.
\end{abstract}
\maketitle
\section{Introduction}
For a compact metric space $(X, d)$, a continuous map $f:X \mapsto X$ and a continuous function $\varphi: X \mapsto \IR$, we continue a program started in \cite{Tho4} to understand the topological pressure of the multifractal decomposition
\[
X = \bigcup_{\alpha \in \IR } X(\varphi, \alpha) \cup \widehat X(\varphi),
\]
where $X(\varphi, \alpha)$ denotes the set of points
\[
X(\varphi, \alpha) = \left \{ x \in X : \lim_{n \ra \infty} \frac{1}{n} \sum_{i = 0}^{n-1} \varphi (f^i (x)) = \alpha \right \}
\]
and $\widehat X(\varphi)$ denotes the set of points for which the Birkhoff average does not exist. In \cite{Tho4}, we showed that $\widehat X(\varphi)$ is either empty or has full topological pressure. In the present work, we turn our attention to the sets $X(\varphi, \alpha)$. Our main result (theorem \ref{theorem1}) is that for any continuous functions $\varphi, \psi: X \mapsto \IR$,
\begin{equation} \label{main}
P_{X(\varphi,\alpha)} (\psi) = \sup \left \{ h_\mu + \int \psi d \mu : \mu \in \mathcal M_f (X) \mbox{ and } \int \varphi d\mu = \alpha \right \},
\end{equation}
where $P_{X(\varphi, \alpha)} (\psi)$ denotes the topological pressure of $\psi$ on $X(\varphi, \alpha)$, defined in \S \ref{pressure}. The motivation for proving multifractal analysis results where pressure is the dimension characteristic is twofold. Firstly, topological pressure is a non-trivial and natural generalisation of topological entropy, which is the standard dynamical dimension characteristic. 
Secondly, understanding the topological pressure of the multifractal decomposition allows us to prove results about the topological entropy of systems related to the original system, for example, suspension flows. 


The class of maps satisfying the specification property includes the time-$1$ map of the geodesic flow of compact connected negative curvature manifolds and certain quasi-hyperbolic toral automorphisms as well as any system which can be modelled by a topologically mixing shift of finite type (see \cite{Tho4} for details).

Formulae similar to (\ref{main}) have a key role in multifractal analysis. For hyperbolic maps and H\"older continous $\varphi$, Barriera and Saussol established our main result for the case $\psi = 0$, i.e. for the topological entropy of $X(\varphi, \alpha)$ and used it to give a new proof of the multifractal analysis in this setting \cite{BS4}. Takens and Verbitskiy proved (\ref{main}) in the case of topological entropy for maps with the specification property \cite{TV}.

Luzia proved our main result for topological pressure when the system is a topologically mixing subshift of finite type and $\varphi, \psi$ are H\"older, and used it to analyse fibred systems \cite{L}.  Our current result generalises and unifies the above mentioned results.

Pfister and Sullivan generalised the result of Takens and Verbitskiy still further, to a setting which applies to $\beta$ shifts \cite{PfS} . We expect that our current method can be extended to their setting. 

Barreira and Saussol proved an analogue of (\ref{main}) for hyperbolic flows when $\psi = 0$ and $\varphi$ is H\"older \cite{BS5}. While we expect (\ref{main}) can be established for flows with specification using our current methods, we consider here the class of suspension flows over maps with specification, and show that (\ref{main}) holds true in this setting.

A large part of our argument is the same as that used by the author in \cite{Tho4}, which was inspired by Takens and Verbitskiy \cite{TV}. We do not give a self-contained proof of this part of the argument but state the key ideas and refer the reader to \cite{Tho4} for the details.

An interesting application of our main result is a `Bowen formula' for the Hausdorff dimension of the level sets of the Birkhoff average for a class of non-uniformly expanding maps of the  interval, which includes the Manneville-Pomeau family of maps.

In \S \ref{sec2}, we take care of our preliminaries. In \S \ref{sec3}, we state and prove our main results. In \S \ref{sec4}, we apply our main result to suspension flows. In \S \ref{sec5}, we use our main result to derive a certain Bowen formula for interval maps.

\section{Preliminaries} \label{sec2}
We give the definitions and fix the notation necessary to give a precise statement of our results, including topological pressure for non-compact sets and the specification property. Let $(X,d)$ be a compact metric space and $f:X \mapsto X$ a continuous map. Let $C(X)$ denote the space of continuous functions from $X$ to $\IR$, and $\varphi, \psi \in C(X)$. Let $S_n \varphi (x) := \sum_{i = 0}^{n-1} \varphi (f^i (x))$ and for $c >0$, let $\mbox{Var}(\varphi, c) := \sup \{ |\varphi(x) - \varphi(y)| : d(x,y)< c\}$. Let $\mathcal{M}_{f} (X)$ denote the space of $f$-invariant probability measures and $\mathcal{M}^e_{f} (X)$ denote those which are ergodic. 
We define the empirical measures
\[
\delta_{x, n} = \frac{1}{n} \sum_{k=0}^{n-1} \delta_{f^k (x)},
\]
where $\delta_x$ is the Dirac measure at $x$.

Given $\epsilon > 0, n \in \IN$ and a point $x \in X$, define the open $(n, \epsilon)$-ball at $x$ by
\[
\mathit{B}_{n}(x, \epsilon) = \{ y \in X : d( f^i(x), f^i(y)) < \epsilon \mbox{ for all } i = 0, \ldots , n-1\}.
\]
Alternatively, let us define a new metric
\[
d_n (x, y) = \max \{ d( f^i(x), f^i(y)) : i = 0, 1, \ldots, n-1 \}.
\]
It is clear that $\mathit{B}_{n}(x, \epsilon)$ is the open ball of radius $\epsilon$ around $x$ in the $d_n$ metric, and that if $n \leq m$ we have $d_n (x, y) \leq d_m (x, y)$ and $\mathit{B}_{m}(x, \epsilon) \subseteq \mathit{B}_{n}(x, \epsilon)$.

Let $Z \subset X$. We say a set $\mc S \subset Z$ is an $(n, \epsilon)$ spanning set for $Z$ if for every $z \in Z$, there exists $x \in \mc S$ with $d_n (x, z) \leq \epsilon$. 
We say a set $\mc R \subset Z$ is an $(n, \epsilon)$ separated set for $Z$ if for every $x, y \in \mc R$, $d_n (x, y) > \epsilon$. 
See \cite{Wa} for the basic properties of spanning sets and seperated sets. 

\subsection{Definition of the topological pressure} \label{pressure}
Let $Z \subset X$ be an arbitrary Borel set, not necessarily compact or invariant. We use the definition of topological pressure as a characteristic of dimension type, due to Pesin and Pitskel. 
We consider finite and countable collections of the form $\Gamma = \{ B_{n_i}(x_i, \epsilon) \}_i$. 
For $\alpha \in \IR$, we define the following quantities:
\[
Q(Z,\alpha, \Gamma, \psi) = \sum_{B_{n_i}(x_i, \epsilon) \in \Gamma} \exp \left(-\alpha n_i +\sup_{x\in B_{n_i}(x_i, \epsilon)} \sum_{k=0}^{n_i -1} \psi(f^{k}(x)) \right),
\]
\[
M(Z, \alpha, \epsilon, N, \psi) = \inf_{\Gamma} Q(Z,\alpha, \Gamma, \psi),
\]
where the infimum is taken over all finite or countable collections of the form $\Gamma = \{ B_{n_i}(x_i, \epsilon) \}_i$ with $x_i \in X$ such that $\Gamma$ covers Z and $n_i \geq N$ for all $i = 1, 2, \ldots$. Define
\[
m(Z, \alpha, \epsilon, \psi) = \lim_{N \rightarrow \infty} M(Z, \alpha,\epsilon, N, \psi).
\]
The existence of the limit is guaranteed since the function $M(Z, \alpha,\epsilon, N)$ does not decrease with N. By standard techniques, we can show the existence of
\begin{displaymath}
P_Z (\psi, \epsilon) := \inf \{ \alpha : m(Z, \alpha, \epsilon, \psi) = 0\} = \sup \{ \alpha :m(Z, \alpha, \epsilon, \psi) = \infty \}.
\end{displaymath}
\begin{definition}
The topological pressure of $\psi$ on $Z$ is given by
\[
P_Z (\psi) = \lim_{\epsilon \ra 0} P_Z (\psi, \epsilon).
\]
\end{definition}
See \cite{Pe} for verification of well-definedness of the quantities $P_Z (\psi, \epsilon)$ and $P_Z (\psi)$. If $Z$ is compact and invariant, our definition agrees with the usual topological pressure as defined in \cite{Wa}. 

\subsection{The specification property}
We are interested in transformations $f$ of the following type: 
\begin{definition} \label{3a}
A continuous map $f: X \mapsto X$ satisfies the specification property if for all $\epsilon > 0$, there exists an integer $m = m(\epsilon )$ such that for any collection $\left \{ I_j = [a_j, b_j ] \subset \IN : j = 1, \ldots, k \right \}$ of finite intervals with $a_{j+1} - b_j \geq m(\epsilon ) \mbox{ for } j = 1, \ldots, k-1 $ and any $x_1, \ldots, x_k$ in $X$, there exists a point $x \in X$ such that
\begin{equation} \label{3a.02}
d(f^{p + a_j}x, f^p x_j) < \epsilon \mbox{ for all } p = 0, \ldots, b_j - a_j \mbox{ and every } j = 1, \ldots, k.
\end{equation}
\end{definition}
The original definition of specification, due to Bowen, was stronger.
\begin{definition} \label{3a.2}
We say $f: X \mapsto X$ satisfies Bowen specification if under the assumptions of definition \ref{3a} and for every $p \geq b_k - a_1 + m(\epsilon)$, there exists a periodic point $x \in X$ of least period $p$ satisfying (\ref{3a.02}).
\end{definition}
One can describe a map $f$ with specification intuititively as follows. For any set of points $x_1, \ldots, x_k$ in $X$, there is an $x \in X$ whose orbit follows given finite pieces of the orbits of the points $x_1, \ldots, x_k$. In this way, one can connect together arbitrary pieces of orbit. If $f$ has Bowen specification, $x$ can be chosen to be a periodic point of any sufficiently large period. 

One can verify that a map with the specification property is topologically mixing. The following converse result holds \cite{Bl}, a recent proof of which is available in \cite{Bu}.
\begin{proposition} [Blokh]
A continuous topologically mixing map of the interval has Bowen specification.
\end{proposition}
Topologically mixing shifts of finite type have specification and factors of systems with specification have specification. 
We give a survey of many interesting examples of maps with the specification property in \cite{Tho4}. 

\subsection{The multifractal spectrum of Birkhoff averages}

For $\alpha \in \IR$, we define 
\[
X(\varphi, \alpha) = \left \{ x \in X : \lim_{n \ra \infty} \frac{1}{n} \sum_{i = 0}^{n-1} \varphi (f^i (x)) = \alpha \right \}.
\]
We define the multifractal spectrum for $\varphi$ to be
\[
\mathcal L_\varphi := \{ \alpha \in \IR : X(\varphi, \alpha) \neq \emptyset \}.
\]
The following lemma (proof included for completeness) is essentially contained in \cite{TV}.
\bl \label{lemma0}
When $f$ has the specification property, $\mathcal L_\varphi$ is a non-empty bounded interval. Furthermore, $\mathcal L_\varphi = \{ \int \varphi d \mu : \mu \in \mathcal{M}_{f} (X) \}$.
\el
\begin{proof}
We first show that $\mathcal L_\varphi = \mathcal I_\varphi$ where $\mathcal I_\varphi = \{ \int \varphi d \mu : \mu \in \mathcal{M}_{f} (X) \}$. By Proposition 21.14 of \cite{De}, when $f$ has the Bowen specification property, every $f-$invariant (not necessarily ergodic) measure has a generic point (i.e. a point $x$ which satisfies $\frac{1}{n}S_n \varphi (x) \ra \int \varphi d \mu$ for all continuous functions $\varphi$). One can verify that this remains true under the specification property. Thus, given $\mu \in \mc M_f (X)$, any choice $x$ of generic point for $\mu$ lies in $X(\varphi, \int \varphi d \mu)$ and so $\mathcal I_\varphi \subseteq \mathcal L_\varphi$. Now take $\alpha \in \mathcal L_\varphi$ and any $x \in X(\varphi,\alpha)$. Let $\mu$ be any weak$^\ast$ limit of the sequence $\delta_{x, n}$. It is a standard result that $\mu$ is invariant, and easy to verify that $\int \varphi d \mu = \alpha$. Thus $\mathcal I_\varphi = \mathcal L_\varphi$.

It is clear that $\mathcal I_\varphi \subseteq [ \inf_{x \in X} \varphi (x), \sup_{x \in X} \varphi (x) ]$ and is non-empty. To show $\mathcal I_\varphi$ is an interval we use the convexity of $\mathcal{M}_{f} (X)$. Assume $\mathcal I_\varphi$ is not a single point. Let $\alpha_1, \alpha_2 \in \mathcal I_\varphi$. Let $\beta \in (\alpha_1, \alpha_2)$. Let $\mu_i$ satisfy $\int \varphi d \mu_i = \alpha_i$ for $i = 1, 2$. Let $t \in (0, 1)$ satisfy $\beta = t \alpha_1 + (1-t) \alpha_2$. One can easily see that $m := t \mu_1 + (1-t) \mu_2$ satisfies $\int \varphi d m = \beta$, and we are done.
\end{proof}
Let $\phi_1, \phi_2 \in C(X)$. We say $\phi_1$ is cohomologous to $\phi_2$
if they differ by a coboundary, i.e. there exists $h \in C(X)$ such that
\[
\phi_1 = \phi_2 + h - h \circ f.
\]
If $\phi_1$ and $\phi_2$ are cohomologous, then $\mathcal L_{\varphi_1}$ equals $\mathcal L_{\varphi_2}$. 
\section{Results} \label{sec3}
\bt \label{theorem1}
Suppose $\varphi, \psi \in C(X, \IR)$ and $\alpha \in \mc L_\varphi$, then
\[
P_{X(\varphi, \alpha)} (\psi) = \sup \left \{ h_\mu + \int \psi d \mu : \mu \in \mathcal M_f (X) \mbox{ and } \int \varphi d\mu = \alpha \right \}.
\]
\et
As a simple corollary, we note that if $\alpha = \int \varphi d m_{\psi}$, where $m_\psi$ is an equilibrium measure for $\psi$ (in the usual sense), then $P_{X(\varphi, \alpha)} (\psi) =P_X(\psi)$.
\subsection{Upper Bound on $P_{X(\varphi,\alpha)} (\psi)$} \label{upper}
We clarify the method of Takens and Verbitskiy. Our proof relies on analysis of the lower capacity pressure of $X(\varphi,\alpha)$, which we define now.
For $Z \subset X$, let
\[
Q_n (Z, \psi, \epsilon) = \inf \left \{ \sum_{x \in S} \exp \left \{ \sum_{k=0}^{n -1} \psi(f^{k}x) \right \}: S \mbox{ is (n, $\epsilon$) spanning set for } Z \right \}, 
\]
\[
P_n (Z, \psi, \epsilon) = \sup \left \{ \sum_{x \in S} \exp \left \{ \sum_{k=0}^{n -1} \psi(f^{k}x) \right \}: S \mbox{ is (n, $\epsilon$) separated set for } Z \right \}. 
\]
We have $Q_n (Z, \psi, \epsilon) \leq P_n (Z, \psi, \epsilon)$. Define
\[
\ul {CP}_{Z} (\psi, \epsilon) = \liminf_{n \ra \infty} \frac{1}{n} \log Q_n (Z, \psi, \epsilon), 
\]
\[
\ul {CP}_{Z} (\psi) = \lim_{\epsilon \ra 0} \ul {CP}_{Z} (\psi, \epsilon).
\]
It is proved in \cite{Pe} that $P_Z(\psi) \leq \ul {CP}_Z(\psi)$. 
\bl \label{uniform}
When $f$ has the specification property, given $\gamma > 0$, there exists $Z \subset X(\varphi, \alpha)$, $t_k \ra \infty$ and $\epsilon_k \ra 0$ such that if $p \in Z$ then
\begin{equation} \label{unif.01}
\left| \frac{1}{m} S_m\varphi (p) - \alpha \right| \leq \epsilon_k \mbox{ for all } m \geq t_k
\end{equation}
and $\ul {CP}_{Z} (\psi) \geq \ul {CP}_{X(\varphi, \alpha)} (\psi) - 4\gamma$.
\el
\bp
Choose $\epsilon > 0$ such that $\ul {CP}_{X(\varphi, \alpha)}(\psi, 2\epsilon) \geq \ul {CP}_{X(\varphi, \alpha)} (\psi) - \gamma$. For $\delta >0$, let
\[
X(\alpha, n, \delta) =\{ x \in X(\varphi, \alpha)  : \left| \frac{1}{m} S_m\varphi (x) - \alpha \right| \leq \delta \mbox{ for all } m \geq n\}.
\] 
We have $X( \varphi, \alpha) = \bigcup_n X(\alpha, n, \delta)$ and $X(\alpha, n, \delta) \subset X(\alpha, n+1, \delta)$, thus $\ul {CP}_{X(\varphi, \alpha)} (\psi, 2\epsilon) = \lim_{n \ra \infty} \ul {CP}_{X(\alpha, n, \delta)} (\psi, 2\epsilon)$. Let $\delta_k \ra 0$ be arbitrary and for each $\delta_k$, pick $M_k \in \IN$ so that
\[
\ul {CP}_{X(\alpha, M_k, \delta_k)} (\psi, 2\epsilon) \geq \ul {CP}_{X(\varphi, \alpha)} (\psi, 2\epsilon) - \gamma.
\]
Write $X_k := X(\alpha, M_k, \delta_k)$. Let $m_k = m(\epsilon /2^k)$ be as in the definition of specification. Now pick a sequence $N_k$ so that $N_{k+1} > \exp \{\sum_{i=1}^k (N_i + m_i)\}$, $N_k > \exp M_{k+1}$, $N_k > \exp m_{k+1}$ and
\[
Q_{N_k}(X_k, \psi, 2\epsilon) >\exp N_k( \ul {CP}_{X(\varphi, \alpha)} (\psi) - 3 \gamma).
\]
Let $t_1 = N_1$ and $t_{k} = t_{k-1} + m_k + N_{k}$ for $k \geq 2$. Note that $t_k/N_k \ra 1$. 

Using the specification property, we define $Z$ to be the set of all points of the form $p := \bigcap_{k=1}^\infty \overline B_{t_{k}} (z_k, \epsilon /2^{k-1})$, where $z_1 \in X_1$, $z_2$ satisfies
\[
d_{N_1}(z_2, z_1) < \epsilon /4 \mbox{ and }d_{N_2} (f^{N_1 +m_2}z_2 , x_2) < \epsilon / 4\] 
for some $x_2 \in X_2$ and $z_k$ satisfies
\[
d_{t_{k-1}}(z_{k-1}, z_k) < \epsilon / 2^k \mbox{ and } d_{N_k} (f^{t_{k-1} +m_k}z_k, x_k) < \epsilon / 2^k
\] 
for some $x_k \in X_k$. We can verify that $\overline B_{t_{k+1}} (z_{k+1}, \epsilon /2^k) \subset \overline B_{t_{k}} (z_k, \epsilon /2^{k-1})$ and so $p$ is well defined.

For $p \in Z$, there is a uniform error term $\epsilon^\prime_k > \delta_k$ which depends on $\delta_k$ and $\mbox{Var}(\varphi, \epsilon/2^k)$ and $t_{k-1}/N_{k-1}$,  so that $| \frac{1}{t_k} S_{t_k} \varphi (p) - \alpha | < \epsilon^\prime_k$. Now let $t_k <n < t_{k+1}$. Suppose $n-t_k+m_k \geq M_{k+1}$. There exists $x \in X_{k+1}$ such that $d_{N_{k+1}}(f^{t_k+m_k} p, x) < \epsilon/2^{k+1}$ and thus
\beq
S_{n} \varphi (p) &\leq& t_k(\alpha + \epsilon^\prime_k) + (n-t_k) (\alpha+ \delta_{k+1} + \mbox{Var}(\varphi, \epsilon/2^{k+1})) +m_{k+1} \| \varphi \|.
\eeq
Suppose $n-t_k < M_{k+1}$. Then
\[
\frac{1}{n} S_{n} \varphi (x) \leq \frac{t_k}{n} (\alpha + \epsilon^\prime_k) + \frac{n-t_k}{n} \| \varphi \| \leq \alpha + \epsilon^\prime_k + \frac{M_{k+1}}{N_k} \| \varphi \|.
\]
Let $\epsilon_{k} = \max\{ \epsilon^\prime_k, \delta_{k+1} + \mbox{Var}(\varphi, \epsilon/2^{k+1})\} + \max \{M_{k+1}/ N_k, m_{k+1}/ N_k\}  \| \varphi\|$ and we have shown that (\ref{unif.01}) holds.

Take a $(t_k, \epsilon)$ spanning set $S_k$ satisfying $\sum_{x \in S_k} \exp S_{t_k} \psi (x) = Q_{t_k} (Z, \psi, \epsilon)$.
It follows that $f^{t_{k-1}+m_k} S_k$ is a $(N_k, \epsilon)$ spanning set for $f^{t_{k-1}+m_k} Z$. Since $\sup \{d_{N_k} (x, z) : x \in X_k, z \in f^{t_{k-1}+m_k} Z \} < \epsilon /2^k$, then $f^{t_{k-1}+m_k} S_k$ is a $(N_k, 2 \epsilon)$ spanning set for $X_k$. Thus
\[
\sum_{x \in S_k} \exp S_{N_k} \psi (f^{t_{k-1}+m_k}x) \geq Q_{N_k} (X_k, \psi, 2\epsilon) > \exp N_k( \ul {CP}_{X(\varphi, \alpha)} (\psi) - 3 \gamma),
\]
and for sufficiently large $k$,
\beq
\sum_{x \in S_k} \exp S_{t_k} \psi (x) &\geq& \exp \{N_k( \ul {CP}_{X(\varphi, \alpha)} (\psi) - 3 \gamma) + (t_{k-1} +m_k) \inf \psi \}\\
&\geq& \exp \{t_k( \ul {CP}_{X(\varphi, \alpha)} (\psi) - 4 \gamma)\}.
\eeq
Taking the $\liminf$ of the sequence $t_k^{-1} \log Q_{t_k} (Z, \psi, \epsilon)$, it follows that
\[
\ul{CP}_Z (\psi, \epsilon) > \ul {CP}_{X(\varphi, \alpha)} (\psi) - 4 \gamma.
\] Since $\epsilon$ was arbitrary, we're done.
\ep
We follow the second half of the proof of the variational principle (Theorem 9.10 of \cite{Wa}). We construct a measure out of $(n, \epsilon)$ separated sets for $Z$ (with a suitable fixed choice of $\epsilon$). 
In contrast, Takens and Verbitskiy construct a measure from $(n, \epsilon_n)$ separated sets with $\epsilon_n \ra 0$. We believe it is not clear in this case how to use the proof of the variational principle to give the desired result. The uniform convergence provided by lemma \ref{uniform} is designed to avoid this. We fix $\gamma >0$ and find $\epsilon >0$ such that $\ul {CP}_Z (\psi, \epsilon) > \ul {CP}_Z (\psi) - \gamma$.
 
Let $S_n$ be a $(n, \epsilon)$ separated set for $Z$ with
\[
\sum_{x \in S_n} \exp S_{n} \psi (x) = P_{n} (Z, \psi, \epsilon),
\]
and write $P_n := P_{n} (Z, \psi, \epsilon)$. 
Let $\sigma_n \in \mc M (X)$ be given by
\[
\sigma_n = \frac{1}{P_n} \sum_{x \in S_k} \exp S_{n} \psi(x) \delta_x
\]
and let 
\[
\mu_n = \frac{1}{n} \sum_{i=0}^{n-1} \sigma_n \circ f^{-i}.
\]
Let $n_j$ be a sequence of numbers so that $\mu_{n_j}$ converges, and let $\mu$ be the limit measure. We have $\mu \in \mc M_f (X)$ and we verify that $\int \varphi d \mu = \alpha$. Let $n \in \IN$ and $k$ be the unique number so $t_k \leq n < t_{k+1}$. Using lemma \ref{uniform}, we have
\beq
\int \varphi d \mu_n &=& \frac{1}{P_n} \frac{1}{n} \sum_{x \in S_k} S_{n} \varphi(x) e^{S_{n} \psi(x)}\\
&\leq& \frac{1}{P_n} \frac{1}{n} \sum_{x \in S_k} n(\alpha +\epsilon_k) e^{S_{n} \psi(x)} \\
&=&\alpha + \epsilon_k,
\eeq
and it follows that $\int \varphi d \mu = \alpha$.

To show that $h_\mu + \int \psi d \mu \geq \liminf_{j \ra \infty} \frac{1}{n_j} \log P_{n_j}$, we recall some key ingredients of the proof of the variational principle, refering the reader to \cite{Wa} for additional notation and details. Let $\xi$ be a partition of $X$ with diameter less than $\epsilon$ and $\mu(\partial \xi)=0$. 
\[
H_{\sigma_n} (\bigvee_{i=1}^{n} f^{-i} \xi) + \int S_{n} \psi d \sigma_n = \log P_n.
\]
Since $\mu(\partial \xi) = 0$, we have for any $k, q \in \IN$,
\[
\lim_{j \ra \infty} H_{\mu_{n_j}} ( \bigvee_{i=0}^{q-1} f^{-i} \xi) = H_{\mu} ( \bigvee_{i=0}^{q-1} f^{-i} \xi).
\]
For a fixed $n$ and $ 1 <q< n$ and $0 \leq j \leq q-1$, we have
\[
\frac{q}{n} \log P_n \leq H_{\mu_{n}} ( \bigvee_{i=0}^{q-1} f^{-i} \xi) +q \int \psi d \mu_n +2\frac{q^2}{n} \log \# \xi. 
\]
Replacing $n$ by $n_j$ and taking $j \ra \infty$, we obtain
\[
q \liminf_{j\ra \infty} \frac{1}{n_j}\log P_{n_j} \leq H_{\mu} ( \bigvee_{i=0}^{q-1} f^{-i} \xi) +q \int \psi d \mu.
\]
Dividing by $q$ and letting $q \ra \infty$, we obtain
\[
\ul {CP}_{Z} (\psi, \epsilon) \leq \liminf_{n \ra \infty} \frac{1}{n} \log P_n \leq h_{\mu} (f, \xi) + \int \psi d \mu \leq h_{\mu}+ \int \psi d \mu.
\]
It follows that 
\[
P_{X(\varphi, \alpha)} (\psi) -5 \gamma \leq \ul {CP}_{X(\varphi, \alpha)} (\psi) -5 \gamma \leq \ul {CP}_{Z} (\psi) - \gamma \leq \ul {CP}_{Z} (\psi, \epsilon) \leq h_\mu + \int \psi d \mu.
\]
Since $\gamma$ was arbitrary, we're done.

\subsection{Lower Bound on $P_{X(\varphi, \alpha)} (\psi)$}
This inequality is harder and the proof is similar to the main theorem of \cite{Tho4}, which we follow closely. The key ingredients are the following two propositions, which respectively generalise the entropy distribution principle \cite{TV} and Katok's formula for measure-theoretic entropy \cite{K}. The first is proved in \cite{Tho4} and the second in \cite{Me}.
\begin{proposition}
\label {theorem3}
Let $f : X \mapsto X$ be a continuous transformation. Let $Z \subseteq X$ be an arbitrary Borel set. Suppose there exists $\epsilon > 0$ and $s \geq 0$ such that one can find a sequence of Borel probability measures $\mu_k$, a constant $K> 0$, and a limit measure $\nu$ of the sequence $\mu_{k}$ satisfying $\nu(Z) > 0$ such that 
\[
\limsup_{k \ra \infty} \mu_{k} (B_n (x, \epsilon)) \leq K \exp \{-ns + \sum_{i=0}^{n-1} \psi(f^i x)\}
\] 
for sufficiently large $n$ and every ball $B_n (x, \epsilon)$ which has non-empty intersection with $Z$. Then $P_Z (\psi, \epsilon) \geq s$.
\end{proposition}
\begin{proposition} \label {theorem4}
Let $(X,d)$ be a compact metric space, $f:X \mapsto X$ be a continuous map and $\mu$ be an ergodic invariant measure. For $\epsilon > 0$, $\gamma \in (0, 1)$ and $\psi \in C(X)$, define 
\[
N^\mu (\psi, \gamma, \epsilon, n) = \inf \left \{ \sum_{x \in S} \exp \left \{ \sum_{i=0}^{n -1} \psi(f^{i}x)\right \} \right \}
\] 
where the infimum is taken over all sets $S$ which $(n, \epsilon)$ span some set $Z$ with $\mu(Z) \geq 1 - \gamma$. We have
\[
h_\mu + \int \psi d \mu  = \lim_{\epsilon \ra 0} \liminf_{n \ra \infty} \frac{1}{n} \log N^\mu (\psi, \gamma, \epsilon, n).
\]
The formula remains true if we replace the $\liminf$ by $\limsup$. 
\end{proposition}
Our strategy is to define a specially chosen family of finite sets $\mc S_k$ using proposition \ref{theorem4}, which will form the building blocks for the construction of a certain fractal $F \subset X_{\varphi, \alpha}$, on which we can define a sequence of measures suitable for an application of proposition \ref{theorem3}.

The first stage of the construction is where our current argument differs from \cite{Tho4}. After this modification, the rest of the construction goes through largely verbatim.
\subsection{Construction of the special sets $\mc S_k$}
Choose a strictly decreasing sequence $\delta_k \ra 0$ and fix an arbitrary $\gamma > 0$. Let us fix $\mu$ satisfying $\int \varphi d \mu = \alpha$ and
\[
h_\mu + \int \psi d\mu \geq \sup \left \{ h_\nu + \int \psi d \nu : \nu \in \mathcal M_f (X) \mbox{ and } \int \varphi d\nu = \alpha \right \} - \gamma.
\]
We cannot assume that $\mu$ is ergodic, so we use the following lemma \cite{Yo}, p.535, to approximate $\mu$ arbitrarily well by convex combinations of ergodic measures.
\bl \label{lemmaA.2}
For each $\delta_k > 0$, there exists $\eta_k \in \mc M_f (X)$ such that $\eta_k = \sum_{i=1}^{j(k)} \lambda_i \eta^k_i$, where $\sum_{i=1}^{j(k)} \lambda_i = 1$ and $\eta^k_i \in \mc M_f^e (X)$, satisfying $ \line \int \varphi d\mu - \int \varphi d \eta_k \line < \delta_k$ and $ h_{\eta_k} > h_\mu - \delta_k$.
\el
Choose a strictly increasing sequence $l_k \ra \infty$ so that each of the sets
\be \label{5}
Y_{k, i} := \left \{ x \in X : \left | \frac{1}{n} S_n \varphi (x) - \int \varphi d \eta^k_i \right | < \delta_k \mbox{ for all } n \geq l_k\right \}
\ee
satisfies $\eta_i^k (Y_{k,i}) > 1 - \gamma$ for every $k \in \IN, i \in \{1, \ldots, j(k)\}$. This is possible by Birkhoff's ergodic theorem. Using proposition \ref{theorem4}, 
we can establish the following lemma (see the corresponding lemma in \cite{Tho4} for details of the proof). Let $\gamma^\prime > 0$.
\begin{lemma} \label {lemma1}
For any sufficiently small $\epsilon > 0$, we can find a sequence $\hat n_k \ra \infty$ with $[\lambda_i \hat n_k] \geq l_k$ and 
finite sets $\mc S_{k,i}$ so that each $\mc S_{k,i}$ is a $([\lambda_i \hat n_k], 5 \epsilon)$ separated set for $Y_{k, i}$ and $M_{k,i} := \sum_{x \in \mc S_{k,i}} \exp \left \{ \sum_{i=0}^{n_k -1} \psi(f^{i}x)\right \}$ satisfies
\[
M_{k,i} \geq \exp\left \{ [\lambda_i \hat n_k] \left(h_{\eta^k_i} + \int \psi d \eta^k_i  - \frac{4}{j(k)} \gamma^\prime \right) \right\}.
\]
Furthermore, the sequence $\hat n_k$ can be chosen so that $\hat n_k \geq  2^{m_{k}}$ where $m_k = m(\epsilon/2^k)$ is as in the definition of specification.
\end{lemma}
We choose $\epsilon$ sufficiently small so that the lemma applies and $\mbox{Var}(\psi, 2 \epsilon) < \gamma$. We fix all the ingredients provided by the lemma. We now use the specification property to define the set $\mc S_k$ as follows. Let $y_i \in \mc S_{k, i}$ and define $x = x(y_1, \ldots, y_{j(k)})$ to be a choice of point which satisfies
\[
d_{[\lambda_i \hat n_k]}(y_l, f^{a_l} x) < \frac{\epsilon}{2^k}
\]
for all $l \in \{1, \ldots, j(k) \}$ where $a_1 = 0$ and $a_{l} = \sum_{i=1}^{l-1} [\lambda_i \hat n_k] + (l-1) m_k$ for $l \in \{2, \ldots, j(k)\}$. Let $\mc S_k$ be the set of all points constructed in this way. Let $n_k = \sum_{i=1}^{j(k)} [\lambda_i \hat n_k] + (j(k)-1) m_k$. Then $n_k$ is the amount of time for which the orbit of points in $\mc S_k$ has been prescribed and we have $n_k/\hat n_k \ra 1$. We can verify that $\mc S_k$ is $(n_k, 4 \epsilon)$ separated and so $\# \mc S_k = \# \mc S_{k,1} \ldots \# S_{k,j(k)}$. Let $M_k := M_{k,1} \ldots M_{k,j(k)}$.

We assume that $\gamma^\prime$ was chosen to be sufficiently small so the following lemma holds.
\bl \label{5.1}
We have

(1) for sufficiently large $k$, $M_k \geq \exp n_k (h_{\mu} + \int \psi d \mu - \gamma)$;

(2) if $x \in \mc S_k$, $|\frac{1}{n_k} S_{n_k} \varphi (x) - \alpha | < \delta_k + \mbox{Var}(\varphi, \epsilon/2^k)+1/k$.
\el
\bp
We have for sufficiently large $k$,
\beq
M_k &\geq& \exp \sum_{i=1}^{j(k)} \{[\lambda_i \hat n_k](h_{\eta^k_i} + \int \psi d \eta^k_i  - 4 j(k)^{-1} \gamma^\prime)\}
\\&\geq& \exp\{(1 - \gamma^\prime)\hat n_k\sum_{i=1}^{j(k)} \lambda_i(h_{\eta_i^k} + \int \psi d \eta_i^k) - 4 \gamma^\prime\}\\
&\geq& \exp (1 - \gamma^\prime)^2 n_k (h_{\eta_k} + \int \psi d \eta_k - 4 \gamma^\prime )\\
&\geq&\exp (1 - \gamma^\prime)^2 n_k (h_{\mu} + \int \psi d \mu - 4 \gamma^\prime -2 \delta_k).
\eeq
Thus if $\gamma^\prime$ is sufficiently small, we have (1).

Suppose $x = x(y_1, \ldots, y_{j(k)}) \in \mc S_k$, then
\beq
|S_{n_k} \varphi (x) - n_k \alpha | &\leq& 
|S_{n_k} \varphi (x) - n_k (\int \varphi d \eta_k - \delta_k)|\\
&\leq&
\sum_{i=1}^{j(k)} |S_{[\lambda_i \hat n_k]} \varphi(f^{a_i} x)- n_k \lambda_i \int \varphi d \eta_i^k|\\
& & +n_k \delta_k + m_k(j(k)-1) \|\varphi\|\\
&\leq& \sum_{i=1}^{j(k)} |S_{[\lambda_i \hat n_k]} \varphi(y_i)- [\lambda_i \hat n_k] \int \varphi d \eta_i^k| + m_k j(k)\|\varphi\| \\& & +n_k \mbox{Var}(\varphi, \epsilon/2_k) +n_k \delta_k \\
&<& \delta_k \sum_{i=1}^{j(k)} [\lambda_i \hat n_k] + m_k j(k)\|\varphi\| + n_k \mbox{Var}(\varphi, \epsilon/2_k) + n_k \delta_k
\eeq
The result follows on dividing through by $n_k$.
\ep

We now construct two intermediate families of finite sets. We follow \cite{Tho4}, to which we refer the reader for the full details. The first such family we denote by $\{ \mathcal C_k \}_{k \in \IN}$ and consists of points which shadow a very large number $N_k$ of points from $\mc S_k$. The second family we denote by $\{ \Tk\}_{k \in \IN}$  and  consist of points which shadow points (taken in order) from ${\mathcal C _1}, {\mathcal C _2}, \ldots, {\mathcal C _k}$. We choose $N_k$ to grow to infinity very quickly, so the ergodic average of a point in $\Tk$ is close to the corresponding point in $\Ck$.
\subsection{Construction of the intermediate sets $\{ \Ck \}_{k \in \IN}$} \label{f}
Let us choose a sequence $N_k$ which increases to $\infty$ sufficiently quickly so that
\begin{equation} \label{f.1}
\limk \frac{n_{k+1} + m_{k+1}}{N_k} = 0,  \limk \frac{N_1 (n_1 + m_1) + \ldots + N_k (n_k + m_k)}{N_{k+1}} = 0.
\end{equation}
We enumerate the points in the sets $\mc S_k$ provided by lemma \ref{lemma1} and write them as follows
\[
\mathcal S_k = \{ x_i^k : i = 1, 2, \ldots, \# \mc S_k \}.
\]
Let us make a choice of $k$ and consider the set of words of length $N_k$ with entries in $\{ 1, 2, \ldots, \# \mc S_k\}$. Each such word $ \ul i = (i_1, \ldots, i_{N_k} )$ represents a point in $\mathcal S_k ^{N_k}$. Using the specification property, 
 we can choose a point $y:= y(i_1, \ldots, i_{N_k} )$ which satisfies
\[
d_{n_k}(x_{i_j}^k, f^{a_j} y) < \frac{\epsilon}{2^k}
\]
for all $j \in \{1, \ldots, N_k \}$ where $a_j = (j-1)(n_k +m_k)$. (i.e. $y$ shadows each of the points $x_{i_j}^k$ in order for length $n_k$ and gap $m_k$.) We define
\[
\Ck = \left \{ y(i_1, \ldots, i_{N_k}) \in X : (i_1, \ldots, i_{N_k}) \in \{1, \ldots, \# \mc S_k \}^{N_k} \right \}.
\]
Let $c_k = N_k n_k + (N_k - 1) m_k$. Then $c_k$ is the amount of time for which the orbit of points in $\Ck$ has been prescribed. It is a corollary of the following lemma that distinct sequences $(i_1, \ldots, i_{N_k} )$ give rise to distinct points in $\Ck$. Thus the cardinality of $\Ck$, which we shall denote by $\# \Ck$, is $\#S_k^{N_k}$.
\bl \label{lemma2}
Let $\ul i$ and $\ul j$ be distinct words in $\{ 1, 2, \ldots M_k \}^{N_k}$. Then $y_1 := y ( \ul i )$ and $y_2 := y ( \ul j )$ are $(c_k, 3 \epsilon)$ separated points (ie. $d_{c_k} ( y_1, y_2 ) > 3 \epsilon$).
\el
\subsubsection{Construction of the intermediate sets $\{ \Tk \}_{k \in \IN}$} \label{g}
We define $\Tk$ inductively. Let $\mathcal T _1 = \mathcal C _1$. We construct $\mc T_{k+1}$ from $\mc T_{k}$ as follows. Let $x \in \Tk$ and $y \in \mc C_{k+1}$. Let $t_1 = c_1$ and $t_{k+1} = t_k + m_{k+1} + c_{k+1}$. Using specification, we can find a point $z := z(x, y)$ which satisfies
\[
d_{t_k}(x, z) < \frac{\epsilon}{2^{k+1}} \mbox{ and } d_{c_{k+1}}(y, f^{t_k +m_{k+1}}z) <\frac{\epsilon}{2^{k+1}} .
\]
Define $\mc T_{k+1} = \{ z(x, y) : x \in \Tk, y \in \mc C_{k+1} \}$.
Note that $t_k$ is the amount of time for which the orbit of points in $\Tk$ has been prescribed. Once again, points constructed in this way are distinct. So we have
\[
\# \Tk = \# \mc C_1 \ldots  \# \mc C_k= \#S_1^{N_1} \ldots \# S_k^{N_k}.
\]
This fact is a corollary of the following straight forward lemma:
\begin{lemma} \label{lemma3}
For every $x \in \Tk$ and distinct $y_1, y_2 \in \mc C_{k+1}$
\[
d_{t_k} (z(x,y_1), z(x, y_2)) < \frac{\epsilon}{2^k} \mbox{ and } d_{t_{k+1}} (z(x,y_1), z(x, y_2)) > 2 \epsilon.
\]
Thus $\Tk$ is a $(t_k, 2 \epsilon)$ separated set. In particular, if $z, z' \in \Tk$, then \[
\overline B_{t_{k}} (z, \frac{\epsilon}{2^k}) \cap \overline B_{t_{k}} (z', \frac{\epsilon}{2^k}) = \emptyset.\]
\end{lemma}
\bl \label{lemma4}
Let $z =z(x,y)\in \Tkk$, then
\[
\overline B_{t_{k+1}} (z, \frac{\epsilon}{2^k}) \subset \overline B_{t_{k}} (x, \frac{\epsilon}{2^{k-1}}).
\]
\el
\subsubsection{Construction of the fractal $F$ and a special sequence of measures $\mu_k$} \label{ss}
Let $F_k = \bigcup_{x \in \Tk} \overline B_{t_k} (x, \frac{\epsilon}{2^{k-1}})$. By lemma \ref{lemma4}, $F_{k+1} \subset F_k$. Since we have a decreasing sequence of connected compact sets, the intersection $F = \bigcap_k F_k$ is non-empty. Further, every point $p \in F$ can be uniquely represented by a sequence $\underline p = (\ul p_1, \ul p_2, \ul p_3,. \ldots )$ where each $\ul p_i = (p_1^i, \ldots, p_{N_i}^i) \in \{ 1, 2, \ldots M_i \}^{N_i}$. Each point in $\Tk$ can be uniquely represented by a finite word $(\ul p_1, \ldots \ul p_k)$. We introduce some useful notation to help us see this. Let $y(\ul p_i) \in \mc C_i$ be defined as in \ref{f}. Let $z_1 (\ul p ) = y(\ul p_1)$ and proceeding inductively, let $z_{i+1} ( \ul p) = z (z_{i} ( \ul p), y(\ul p_{i+1})) \in \mc T _{i+1}$ be defined as in \ref{g}. We can also write $z_i (\ul p)$ as $z(\ul p_1, \ldots, \ul p_i)$.  Then define $p := \pi \ul p$ by
\[
p = \bigcap_{i\in \IN} \overline B_{t_{i}} (z_{i} ( \ul p), \frac{\epsilon}{2^{i-1}}).
\]
It is clear from our construction that we can uniquely represent every point in $F$ in this way.
\bl \label{dist}
Given $z = z(\ul p_1, \ldots, \ul p_k) \in \Tk$, we have for all $i \in \{1, \ldots, k\}$ and all $l \in \{1, \ldots, N_i\}$,
\[
d_{n_i}(x^i_{p^i_{l}}, f^{t_{i-1} + m_{i-1} + (l-1)(m_i + n_i)}z) < 2 \epsilon.
\]
\el
We now define the measures on $F$ which yield the required estimates for the Pressure Distribution Principle. For each $z \in \Tk$, we associate a number $\mc L (z) \in (0, \infty)$. Using these mumbers as weights, we define, for each $k$, an atomic measure centred on $\Tk$. Precisely, if $z = z(\ul p_1, \ldots \ul p_k)$, we define
\[
\mc L(z) := \mc L(\ul p_1) \ldots \mc L(\ul p_k),
\]
where if $\ul p_i = (p^i_{1}, \ldots, p^i_{N_i})\in \{1, \ldots, \# \mc S_i \}^{N_i}$, then
\[
\mc L(\ul p_i) := \prod_{l=1}^{N_i} \exp S_{n_i} \psi (x^i_{p^i_l}). 
\]
We define
\[
\nu_k := \sum_{z \in \Tk}\delta_z \mc L(z). 
\]
We normalise $\nu_k$ to obtain a sequence of probability measures $\mu_k$. More precisely, we let $\mu_k := \frac{1}{\kappa_k} \nu_k$, where $\kappa_k$ is the normalising constant
\[
\kappa_k := \sum_{z \in \Tk} \mc L_k (z).
\]
\bl
$\kappa_k = M_1^{N_1} \ldots M_k^{N_k}$.
\el
\bl \label{lemma5}
Suppose $\nu$ is a limit measure of the sequence of probability measures $\mu_k$. 
Then $\nu (F) = 1$.
\el
In fact, the measures $\mu_k$ converge. However, by using the generalised pressure distribution principle, we do not need to use this fact and so we omit the proof (which goes like lemma 5.4 of \cite{TV}). The proof of the following lemma is similar to lemma 5.3 of \cite{TV} or the corresponding lemma of \cite{Tho4}, and relies on (2) of lemma \ref{5.1}. 
\bl \label{lemma6}
For any $p \in F$, the sequence $\lim_{k \ra \infty} \frac{1}{t_k} \sum_{i = 0}^{t_k-1} \varphi (f^i (p)) = \alpha$. Thus $F \subset X(\varphi, \alpha)$.
\el 
For an affirmative answer to theorem \ref{theorem1}, we give a sequence of lemmas which will allow us to apply the generalised pressure distribution principle. The proofs are the same as the corresponding lemmas from \cite{Tho4}, with minor modifications coming from the changed definition of $\mc S_k$ and lemma \ref{5.1}.

Let $\mc B := B_n (q, \epsilon /2)$ be an arbitrary ball which intersects $F$. Let $k$ be the unique number which satisfies $t_k \leq n < t_{k+1}$. Let $j \in \{0, \ldots, N_{k+1} -1 \}$ be the unique number so
\[ 
t_k + (n_{k+1} + m_{k+1})j \leq n < t_k + (n_{k+1} + m_{k+1})(j+1).
\]
We assume that $j \geq 1$ and leave the details of the simpler case $j=0$ to the reader. The following lemma reflects the fact that the number of points in $\mc B \cap \mc T_{k+1}$ is restricted since $\mc T_k$ is $(t_k, 2 \epsilon)$ separated and $\mc S_{k+1}$ is $(n_{k+1}, 4 \epsilon)$ separated. 
\bl \label{lemma6.11}
Suppose $\mu_{k+1} (\mc B) > 0$, then there exists (a unique choice of) $x \in \Tk$ and $i_1, \ldots, i_j \in \{1, \ldots, \# \mc S_{k+1} \}$ satisfying
\[
\nu_{k+1} (\mc B) \leq \mc L (x) \prod_{l=1}^j \exp S_{n_{k+1}} \psi(x^{k+1}_{i_l}) M_{k+1}^{N_{k+1}-j}.
\]
\el
The following lemma is a consequence of lemma \ref{dist}.
\bl \label{lkx}
Let $x \in \Tk$ and $i_1, \ldots, i_j$ be as before. Then
\beq
\mc L(x) \prod_{l=1}^j \exp S_{n_{k+1}} \psi(x^{k+1}_{i_l}) \leq \exp \{S_n \psi(q) &+& 2n \mbox{Var}(\psi, 2 \epsilon)\\&+&\| \psi \| (\sum_{i=1}^k N_i m_i  +j m_{k+1} ) \}.
\eeq
\el
The following lemma 
reflects the restriction on the number of points that can be contained in $\mc B \cap \mc T_{k+p}$.

\bl \label{mk+p}
For any $p \geq 1$, suppose $\mu_{k+p} (\mc B) > 0$. Let $x \in \Tk$ and $i_1, \ldots, i_j$ be as before. 
We have
\[
\nu_{k+p} (\mc B) \leq \mc L (x) \prod_{l=1}^j \exp S_{n_{k+1}} \psi (x^{k+1}_{i_l}) M_{k+1}^{N_{k+1}-j} M_{k+2}^{N_{k+2}} \ldots M_{k+p}^{N_{k+p}}.
\]
\el
\bl \label{8.4}
\[\mu_{k+p} (\mc B) \leq \frac{1}{\kappa_k M_{k+1}^j}\exp \left \{S_n \psi(q) + 2n Var (\psi, 2 \epsilon) +\| \psi \| (\sum_{i=1}^k N_i m_i  +j m_{k+1} ) \right \}. \]
\el


Let $C := h_\mu + \int \varphi d \mu$. The following lemma is implied by lemma \ref{5.1}.
\begin{lemma} \label{lemma8}
For sufficiently large $n$, $\kappa_k M_{k+1}^j \geq \exp((C - 2 \gamma)n)$
\end{lemma}
Combining the previous two lemmas gives us
\bl 
For sufficiently large $n$,
\[
\limsup_{l \ra \infty} \mu_{l} (B_n (q, \frac{\epsilon}{2})) \leq \exp \{-n(C - 2Var (\psi, 2 \epsilon) - 3\gamma) + \sum_{i=0}^{n-1} \psi(f^i q)\}.
\] 
\el
Applying the Generalised Pressure Distribution Principle, 
we have
\[
P_F (\psi, \epsilon) \geq C - 2\mbox{Var} (\psi, 2 \epsilon) - 3\gamma.
\] 
Recall that $\epsilon$ was chosen sufficiently small so $\mbox{Var}(\psi, 2 \epsilon) < \gamma$. It follows that
\[
P_{X(\varphi, \alpha)}(\psi, \epsilon) \geq P_F (\psi, \epsilon) \geq C - 5 \gamma.
\]
Since $\gamma$ and $\epsilon$ were arbitrary, the proof of theorem \ref{theorem1} is complete.
\section{Application to Suspension Flows} \label{sec4}
We apply our main result to suspension flows. Let $f : X \mapsto X$ be a homeomorphism of a compact metric space $(X, d)$. We consider a continuous roof function $\rho : X \mapsto (0, \infty)$. 
We define the suspension space to be
\[
X_\rho = \{ (x, s) \in X \times \IR : 0 \leq s \leq \rho(x) \},
\]
where $(x, \rho(x))$ is identified with $(f(x), 0)$ for all $x$. 
We define the flow $ \Psi = \{ g_t \}$ on $X_\rho$ locally by $g_t (x, s) = (x, s+t)$. To a function $\Phi :X_\rho \mapsto \IR$, we associate the function $\varphi: X \mapsto \IR$ by $\varphi(x) = \int_0^{\rho(x)} \Phi (x, t) dt$. Since the roof function is continuous, when $\Phi$ is continuous, so is $\varphi$. We have (see \cite{Tho4})
\[
\liminf_{T \ra \infty} \frac{1}{T} \int_0^T \Phi(g_t (x,s)) d t = \liminf_{n \ra \infty} \frac{S_n \varphi (x)}{S_n \rho (x)},
\]
\[
\limsup_{T \ra \infty} \frac{1}{T} \int_0^T \Phi(g_t (x,s)) d t = \limsup_{n \ra \infty} \frac{S_n \varphi (x)}{S_n \rho (x)}.
\]
We consider
\beq
X_\rho(\Phi, \alpha) &:=& \{ (x,s) \in \xr : \lim_{T \ra \infty} \frac{1}{T} \int_0^T \Phi(g_t (x,s)) d t = \alpha \}\\
&=& \{ (x, s) : \limn \frac{S_n \varphi (x)}{S_n \rho (x)} = \alpha, 0 \leq s < \rho (x) \}.
\eeq
For $\mu \in \mathcal{M}_{f} (X)$, we define the measure $\mu_\rho$ by
\[
\int_{X_\rho} \Phi d \mu_\rho = \int_X \varphi d \mu / \int \rho d \mu
\]
for all $\Phi \in C (X_\rho)$, where $\varphi$ is defined as above. We have $\Psi$-invariance of $\mu_\rho$ (ie. $\mu (g_t ^{-1} A) = \mu (A)$ for all $t \geq 0$ and measurable sets $A$). The map $\mc R : \mathcal{M}_{f} (X) \mapsto \mathcal{M}_{\Psi} (X_\rho)$ given by $\mu \mapsto \mu_\rho$ is a bijection. It is verified in \cite{PP} that $h_{\mu_\rho} = h_\mu / \int \rho d \mu$ and hence, 
\[
\htop (\Psi) = \sup \{ h_\mu : \mu \in \mathcal{M}_{\Psi} (X_\rho) \} = \sup \left \{\frac{ h_\mu }{ \int \rho d \mu} : \mu \in \mathcal{M}_{f} (X) \right \},
\]
where $\htop(\Psi)$ is the topological entropy of the flow. 
We use the notation  $h_{top} ( Z, \Psi )$ for topological entropy of a non-compact subset $Z \subset X_\rho$ with respect to $\Psi$ (defined in \cite{Tho4}).

\begin{theorem} \label {ratio}
Let $(X,d)$ be a compact metric space and $f:X \mapsto X$ be a continuous map with specification. Let $\varphi, \psi \in C(X)$ and $\rho : X \mapsto (0, \infty)$ be continuous. Let $X (\varphi, \rho, \alpha):=\left \{ x \in X : \limn \frac{S_n \varphi (x)}{S_n \rho (x)} = \alpha \right \}$. 
For $\alpha$ such that $X (\varphi, \rho, \alpha) \neq \emptyset$, we have 
\[
P_{X (\varphi, \rho, \alpha)} (\psi) = \sup \left \{ h_\mu + \int \psi d \mu : \mu \in \mathcal M_f (X) \mbox{ and } \frac{\int \varphi d\mu}{\int \rho d \mu} = \alpha \right \}.
\]
\end{theorem}
\bp
We require only a small modification to the proof of theorem \ref{theorem1}. We modify lemma \ref{lemmaA.2} so $\eta_k$ satisfies $|\int \varphi d\mu / \int \rho d \mu - \int \varphi d\eta_k / \int \rho d \eta_k| < \delta_k$ and replace the family of sets defined at (\ref{5}) by the following:
\[
Y_{k, i} := \left \{ x \in X : \left | \frac{S_n \varphi (x)}{S_n \rho (x)} - \frac{\int \varphi d \eta^k_i}{\int \rho d \eta^k_i}  \right | < \delta_k \mbox{ for all } n \geq l_k\right \}
\]
chosen to satisfy $\eta^k_i (Y_{k,i}) > 1 - \gamma$ for every $k$. This is possible by the ratio ergodic theorem. The rest of the proof requires only superficial modifications.
\ep
\bt \label{lower}
Let $(X,d)$ be a compact metric space and $f:X \mapsto X$ be a homeomorphism with the specification property. Let $\rho : X \mapsto (0, \infty)$ be continuous. Let $(\xr, \Psi)$ be the corresponding suspension flow over $X$. Let $\Phi: \xr \mapsto \IR$ be continuous. We have 
\[
\htop( X_\rho(\Phi, \alpha), \Psi) = \sup \left \{ h_\mu : \mu \in \mathcal M_\Psi (X_\rho) \mbox{ and } \int \Phi d\mu = \alpha \right \}.
\]
\et
\bp
Let $Z \subset X$ be arbitrary and $Z_\rho := \{ (x, s) : x \in Z, 0\leq s < \rho(x) \}$. In \cite{Tho4}, we proved that if $\beta$ is the unique solution  to the equation $P_Z ( - t \rho) = 0$, then $\htop (Z_\rho, \Psi) \geq \beta$. Thus, if $h$ be the unique positive real number which satisfies $P_{X(\varphi, \rho, \alpha)} (-h \rho) = 0$, then $\htop( X_\rho(\Phi, \alpha), \Psi) \geq h$. By theorem \ref{ratio},
\[
\sup \left \{ h_\mu - h \int \rho d \mu : \mu \in \mathcal M_f (X) \mbox{ and } \frac{\int \varphi d\mu}{\int \rho d \mu} = \alpha \right \} = 0.
\]
Thus, if $\mu \in \mathcal M_f (X)$ satisfies $\frac{\int \varphi d\mu}{\int \rho d \mu} = \alpha$, then $h \geq \frac{h_\mu}{\int \rho d \mu}$ and
\beq
h &\geq& \sup \left \{ \frac{h_\mu}{\int \rho d \mu} : \mu \in \mc M_f (X), \frac{\int \varphi d\mu}{\int \rho d \mu} = \alpha \right \} \\
&=& \sup \left \{ h_\mu : \mu \in \mathcal M_\Psi (X_\rho) \mbox{ and } \int \Phi d\mu = \alpha \right \}.
\eeq

For the opposite inequality, we note that $\htop(Z, \Psi) \leq \ul {CP}_{Z} (0)$, where $\ul {CP}_{Z} (0)$ is defined with respect to the time-1 map of $\Psi$. Given $\gamma>0$, we can adapt lemma \ref{uniform} to find a set $Z \subset X_\rho$, $t_k \ra \infty$ and $\epsilon_k \ra 0$ such that for $(x,s) \in X_\rho$, we have
\[
\left|  \frac{1}{T} \int_0^T \Phi(g_t (x,s)) d t- \alpha \right| \leq \epsilon_k \mbox{ for all } T \geq t_k
\]
and $\ul {CP}_{Z} (0) \geq \ul {CP}_{X(\Phi, \alpha)} (0) - 4\gamma$. 
We repeat the argument of \ref{upper} to construct a suitable measure out of $(n, \epsilon)$ spanning sets for the time-$1$ map of the flow which satisfies $\int \Phi d \mu = \alpha$ and $\ul {CP}_{Z} (0) - \gamma \leq h_\mu$. We obtain
\[
\htop( X_\rho(\Phi, \alpha), \Psi) \leq \sup \left \{ h_\mu : \mu \in \mathcal M_\Psi (X_\rho) \mbox{ and } \int \Phi d\mu = \alpha \right \}.\]
\ep
As a simple corollary, we note that if $\alpha = \int \Phi d m$, where $m$ is a measure of maximal entropy for the flow, then $\htop( X_\rho(\Phi, \alpha), \Psi) = \htop (\Phi)$.
\section{A Bowen formula for Hausdorff dimension of level sets of the Birkhoff average for certain interval maps} \label{sec5}
The following application was described to the author by Thomas Jordan. If $f$ is a $C^{1+\alpha} $, uniformly expanding Markov map of the interval and $\varphi: [0,1] \mapsto \IR$, then it was shown by Olsen \cite{Ol} that
\begin{equation} \label{dimh}
\mbox{dim}_H (X(\varphi, \alpha)) = \sup \left \{ \frac{h_\mu}{\int \log f^\prime d \mu} : \int \varphi d \mu = \alpha \right \}.
\end{equation}
In \cite{JP}, the authors consider piecewise $C^1$ Markov maps of the interval with a finite number of parabolic fixed points $x_i$ such that $f(x_i)=x_i$, $f^\prime(x_i) =1$ and $f^\prime (x) > 1$ for $x \in [0,1] \setminus \bigcup_i x_i$. They show that (\ref{dimh}) holds for $\alpha \in \mc L_\varphi \setminus [\min_i\{ \varphi (x_i)\}, \max_i \{\varphi (x_i)\} ]$. 
Simple examples in this category are provided by the Manneville-Pomeau family of maps $f_t (x) = x^t + x^{1+t} (\mbox{mod} 1)$ (where $t > 0$ is a fixed parameter), which have a single parabolic fixed point at $0$. 
Henceforth, we let $\psi = \log f^\prime$. Note that since $\psi$ is non-negative, $s \mapsto P_{X(\varphi, \alpha)} (-s \psi)$ is decreasing (although possibly not strictly decreasing). 
\bt \label{bowen}
Suppose $s \mapsto P_{X(\varphi, \alpha)} (-s \psi)$ has a unique zero $d$ and (\ref{dimh}) holds true. Then $ d= \mbox{dim}_H (X(\varphi, \alpha))$.
\et
\bp
By (\ref{dimh}), if $\mu \in \mc M_f (X)$ and $\int \varphi d \mu = \alpha$, then 
\[
h_\mu - \mbox{dim}_H (X(\varphi, \alpha)) \int \psi d \mu \leq 0.
\] 
By theorem \ref{theorem1}, $P_{X(\varphi, \alpha)} (- \mbox{dim}_H (X(\varphi, \alpha)) \psi) \leq 0$. Thus $\mbox{dim}_H (X(\varphi, \alpha)) \geq d$.

Now suppose $\mbox{dim}_H (X(\varphi, \alpha)) < d$. Since $s \mapsto P_{X(\varphi, \alpha)} (-s \psi)$ is decreasing and has a unique zero, $P_{X(\varphi, \alpha)} (- \mbox{dim}_H (X(\varphi, \alpha)) \psi) > 0$. By theorem \ref{theorem1}, there exists $\mu$ with $\int \varphi d \mu = \alpha$ and $h_\mu - \mbox{dim}_H (X(\varphi, \alpha)) \int \psi d \mu > 0$. This implies that $\mbox{dim}_H (X(\varphi, \alpha)) < h_\mu / \int \psi d \mu$, which contradicts (\ref{dimh}).
\ep
We remark that by a slight modification to the proof, a more general statement is that if (\ref{dimh}) holds and $d = \inf \{s : P_{X(\varphi, \alpha)} (-s \psi) = 0 \}$, then $ d= \mbox{dim}_H (X(\varphi, \alpha))$.

We comment on the hypotheses of theorem \ref{bowen}. If there exists $\mu$ with $\int \varphi d \mu = \alpha$ and $\int \psi d \mu > 0$, then $s \mapsto P_{X(\varphi, \alpha)} (-s \psi)$ is strictly decreasing. Now suppose $\varphi = \psi = \log f^\prime$. In the case of the Manneville-Pomeau family of maps, the only measure with $\int \psi d \mu = 0$ is the Dirac measure supported at $0$, and so $s \mapsto P_{X(\varphi, \alpha)} (-s \psi)$ is decreasing for $\alpha \in \mc L_\varphi \setminus \{ 0 \}$. By \cite{JP}, (\ref{dimh}) holds true for the same set of values and thus theorem \ref{bowen} applies. 
We remark that for $\alpha = 0$, $P_{X(\log f ^\prime, 0)} (-s \psi) = 0$ for all $s \in \IR$.

\section*{Acknowledgements}
This work constitutes part of my PhD, which is supported by the EPSRC. I would like to thank my supervisors Mark Pollicott and Peter Walters for many useful discussions and reading draft versions of this work, for which I am most grateful. I would also like to thank Thomas Jordan for suggesting the final section of this work.

\bibliographystyle{amsplain}
\bibliography{vpforpressure}

\end{document}